\documentclass[12pt]{article}%
   
\usepackage{amsmath,enumerate}
\usepackage{amsfonts}
\usepackage{amssymb}

\newtheorem{theorem}{Theorem}[section]

\newtheorem{claim}[theorem]{Claim}

\newtheorem{corollary}[theorem]{Corollary}

\newtheorem{lemma}[theorem]{Lemma}

\newtheorem{proposition}[theorem]{Proposition}

\newtheorem{question}[theorem]{Question}

\newenvironment{proof}[1][Proof]{\noindent\textbf{#1.} }
{\hfill \ \rule{0.5em}{0.5em}}

\title{Generalized sum-free sets and cycle saturated regular graphs}
\author{David Davini\thanks{Department of Mathematics, University of California, Los Angeles, U.S.A.
\texttt{daviddavini@ucla.edu}}
\and 
Craig Timmons\thanks{Department of Mathematics and Statistics,  California State University Sacramento,  U.S.A. \texttt{craig.timmons@csus.edu}. Research is supported in part by Simons Foundation Grant \#359419}}
\date{\today}

\begin{document}

\maketitle

\begin{abstract}
Gerbner, Patk\'{o}s, Tuza, and Vizer recently 
initiated the study of $F$-saturated regular graphs.
One of the essential problems in 
this line of research is
determining when such a graph exists.     
Using generalized 
sum-free sets
we prove that for any odd integer $k \geq 5$, there is an $n$-vertex regular $C_k$-saturated 
graph for all $n \geq n_k$.  
Our proof is based on constructing
a special type of sum-free set in 
$\mathbb{Z}_n$.   
We prove that for all even $\ell \geq 4$ and integers $n > 12 \ell^2 + 36 \ell + 24$, there is a symmetric complete $( \ell , 1)$-sum-free set in $\mathbb{Z}_n$.  
We pose the problem of finding the minimum size of such a set, and present some examples found by a computer 
search.

\end{abstract}


\section{Introduction}

A graph $G$ is \emph{$F$-free} if $G$ does not contain a subgraph that is isomorphic to $F$.
The graph $F$ is often called the \emph{forbidden subgraph}.
An important class of $F$-free graphs are those 
that are maximal with respect to adding edges.  We say that 
a graph $G$ is \emph{$F$-saturated} if $G$ is $F$-free and 
adding any missing edge to $G$ creates a subgraph that is isomorphic to $F$.  One of the most studied 
problems in graph saturation 
is determining the minimum 
number of edges in an $n$-vertex $F$-saturated graph.
This minimum is called the \emph{saturation number of
$F$}, denoted $\textup{sat}(n,F)$.  
K\'{a}szonyi and Tuza \cite{kt} proved 
that $\textup{sat}(n, F) = O(n)$ for any graph $F$ having at 
least one edge.  A famous result of Erd\H{o}s, Hajnal, Moon \cite{ehm} and Zykov \cite{zykov} 
gives an exact formula for $\textup{sat}(n, K_r)$ for 
$2 \leq r \leq n$, and shows that 
the join of a clique on $r-2$ vertices with 
an independent set on $n -r + 2$ vertices
is the unique extremal graph.
There has been much research on determining saturation numbers of other graphs. 
The surveys of Faudree, Faudree, Schmitt \cite{ffs} and Pikhurko \cite{p} contain
a wealth of information on saturation in graphs and hypergraphs.

A common theme in saturation is to add degree 
restrictions on the family of $F$-saturated graphs; see \cite{aehk, DH, fs, hs} to name a few.
A recent variant of $\textup{sat}(n , F)$ was introduced by 
Gerbner, Patk\'{o}s, Tuza, and Vizer \cite{Gerbner-et-al}.  
Given a graph $F$ and a positive integer $n$, 
one can ask if there exists 
an $n$-vertex $F$-saturated regular graph.  
When such a graph exists, define
\[
\textup{rsat}(n, F)
\]
to be the minimum number of edges in such a graph.  
This can be viewed as a
regular version of saturation numbers much like 
the recently studied regular Tur\'{a}n numbers
\cite{cambie, caro-tuza, Gerbner2, tait}.  

Gerbner 
et al.\ \cite{Gerbner-et-al} proved that $\textup{rsat}(n, K_3)$ exists for all 
sufficiently large $n$.  This can also be proved using 
results of Haviv and Levy \cite{hl} on symmetric complete sum-free sets 
in $\mathbb{Z}_n$.  Among several other interesting theorems, 
Gerbner et al.\ obtained some partial results 
for $K_4$.  The second author \cite{t} proved 
that $\textup{rsat}(n, K_4)$ and $\textup{rsat}(n, K_5)$ exist
for all sufficiently large $n$.  However, the general
problem for larger cliques remains open.  
In \cite{t} it was proved that for all $r \geq 6$, 
$\textup{sat}(n, K_r)$ exists for infinitely many $n$.

Motivated by the approach using symmetric complete sum-free sets, the second author obtained the following theorem on odd cycles.

\begin{theorem}[\cite{t}]\label{old theorem}For all positive integers 
$\alpha$ and $k$, there is a regular $C_{2 \alpha + 3}$-saturated 
graph with $2 \alpha( \alpha + 4) k + 2 \alpha + 5$ vertices.
\end{theorem}

This theorem implies that for any 
odd integer $2 \alpha + 3 \geq 5$,  
there are infinitely many $n$ for which an $n$-vertex $C_{ 2 \alpha + 3}$-saturated regular graph exists.  
However, for fixed $\alpha$ Theorem 
\ref{old theorem} requires $n \equiv 
2 \alpha + 5 \pmod{2 \alpha ( \alpha + 4)} $
which is only one of the possible
$2 \alpha ( \alpha + 4)$
residue classes.  Our first theorem removes 
this restriction.

\begin{theorem}\label{main theorem}
Let $\ell \geq 4$ be an even integer.  For all 
$n > 12 \ell^2 + 36 \ell + 24$, there is an $n$-vertex 
$C_{ \ell + 1}$-saturated regular graph.  
\end{theorem}

Our proof of Theorem \ref{first theorem} is an explicit construction and gives the upper bound
\[
\textup{rsat}(n , C_{ \ell + 1} ) \leq \frac{n^2}{ 2 ( \ell + 1) } - n
\]
for all even $\ell \geq 4$ and odd $n > 12 \ell^2 + 36 \ell + 24$.  For even $n \geq \ell + 1$, a balanced complete bipartite graph is regular and $C_{\ell + 1}$-saturated.  Thus, we have a quadratic upper bound on $\textup{rsat}(n, C_{ \ell +1})$ for all $n > 12 \ell^2 + 36 \ell + 24$ for
even $\ell \geq 4$.
The lower bound $\textup{rsat}(n , C_{ \ell + 1}) = \Omega ( n^{1 + 1/ \ell } )$, when this regular saturation number exists,
is proved in \cite{Gerbner-et-al}.  As discussed in \cite{t}, the construction of
 Haviv and Levy implies 
 $\textup{rsat}(n , C_3) = O  ( n^{3/2})$, and 
 so $\textup{rsat}(n, C_3) = \Theta( n^{3/2} )$.  Determining 
if $\textup{rsat}(n , C_{ \ell + 1} )$ is subquadratic ($\ell \geq 4$ even) is an interesting problem.  An
answer to the following question is a possible first step towards a solution.  

\begin{question}\label{question 1}
Is it true that $\textup{rsat}(n, C_5) = o(n^2)$?  
\end{question}

Our approach follows that of \cite{t} where the 
idea is to use Cayley graphs of $(\ell   , 1 )$-sum-free sets.  
These sets have been studied in additive combinatorics 
\cite{Bajnok, Bajnok-Matzke, HP}, and are a 
generalization of classical sum-free sets (see
the survey of Tao \cite{tv}).  
Let $S$ be a subset of an abelian group $\Gamma$.
For a positive integer
$\ell$, the set $\ell S$ is the \emph{$\ell$-fold sumset of $S$}: \[
\ell S = \{ s_1 + s_2 + \dots + s_{ \ell } : s_i \in S \}.
\]
The set $S$ is \emph{$(k , \ell)$-sum-free} if 
$( k S) \cap ( \ell S) = \emptyset$.  
A $(k , \ell)$-sum-free set is \emph{complete}
if $k S \cup \ell S$ forms a partition of $\Gamma$.  
Finally, $S$ is \emph{symmetric} if $s \in S$ implies $-s \in S$.
The connection between 
$( \ell , 1)$-sum-free sets and regular $C_{ \ell + 1}$-saturated
graphs is discussed in detail in \cite{t}.  Roughly speaking,
one can use a symmetric complete $( \ell , 1)$-sum-free 
set to construct a Cayley graph that will be $C_{ \ell + 1}$-saturated.
This will be made more precise later, but for now, we state 
our main result on $( \ell , 1)$-sum-free sets.  
It is the key ingredient in the proof of Theorem \ref{main theorem}.

\begin{theorem}\label{first theorem}
Let $\ell \geq 4$ be even, $t \geq 1$, $\gamma$ be the unique integer for which
\[
  2\ell + 2 \geq \ell + 3 - 4 ( t + 2) + \gamma ( 2 \ell + 2 ) \geq 
 1,
\]
$j = 2 \gamma + 1$, and $r = \ell + 3 - 4 ( t + 2) + \gamma ( 2 \ell + 2)$.  
If $k $ is any integer with $k \geq 4 t + 2 \ell + 2 | \gamma |+2$, then there is a symmetric complete
$(\ell , 1)$-sum-free set $S \subseteq \mathbb{Z}_{ ( 2 \ell + 2)k + r }$.    
\end{theorem}

The following corollary is simpler to state and is
a byproduct of our proof of Theorem \ref{main theorem} using 
Theorem \ref{first theorem}.  

\begin{corollary}\label{corollary}
Let $\ell \geq 4$ be even.  If  $n > 12 \ell^2 + 36 \ell + 24$,
then there is a symmetric complete $(\ell , 1)$-sum-free set $S \subseteq \mathbb{Z}_n$.  
\end{corollary}

Finding the smallest size of a symmetric complete $( \ell , 1)$-sum-free set in $\mathbb{Z}_n$ 
appears to be a challenging problem.  Write $\psi_{ \ell }(n)$ for this minimum.
Note that $\ell$ must be even for this function to be well-defined.  Indeed,
if $S$ is a non-empty symmetric set in $\mathbb{Z}_n$ and $s \in S$, then 
\[
\underbrace{  s - s  +  s - s + \dots + s - s }_{ 2 \ell ' ~ \mbox{terms in all} } + s \equiv s ( \textup{mod}~n).
\]
This shows that we cannot have 
$( ( 2 \ell' +1) S ) \cap S = \emptyset$ and therefore,
 $\ell$ must be even.  Furthermore, even when restricting $\ell$ to be even
 there may be integers $n \geq 1$ for which $\mathbb{Z}_n$ 
 does not contain a symmetric complete $(\ell,1)$-sum-free set.  For these $n$, $\psi_{ \ell }(n)$ is 
 undefined.  
 
The case $\ell = 2$ corresponds to classical sum-free sets and
$\psi_{ 2 }(n) = \Theta ( n^{1/2} )$ where the upper bound 
is due to Haviv and Levy \cite{hl}.  If $S \subset \mathbb{Z}_n$ is 
a complete $( \ell , 1)$-sum-free set, then 
$\mathbb{Z}_n$ is the disjoint union of $\ell S$ and $S$ so
\[
n \leq \binom{ |S| + \ell - 1}{ \ell } + |S|.
\]
This inequality gives the lower bound $  \psi_{ \ell } (n) =\Omega ( n^{1 / \ell } )$.
The upper bound $\psi_{ \ell }(n) = O(n)$ follows from Corollary \ref{corollary}.

The value of $\psi_4 (n)$ was computed for small values of $n$. 
Our results from $41$ to $80$ are summarized in the table below. For $81 \leq n \leq 140$, our program found that $\psi_4 (n) = 8$ with the exception 
of $n \in \{ 113, , 116, 117, 125 \}$, where $\psi_4 (n) = 10$ for these $n$.

\begin{center}

\begin{tabular}{l|l|l||l|l|l}
$n$ & $\psi_4(n)$ & Example                   & $n$ & $\psi_4(n)$ & Example                    \\ \hline
41  & 6           & 1,5,11,30,36,40           & 61  & 8           & 1,3,5,22,39,56,58,60       \\
42  & 6           & 1,5,18,24,37,41           & 62  & 7           & 1,5,18,31,44,57,61         \\
43  & 6           & 1,6,8,35,37,42            & 63  & 6           & 1,24,28,35,39,62           \\
44  & 6           & 1,7,18,26,37,43           & 64  & 8           & 1,5,9,30,34,55,59,63       \\
45  & 6           & 1,6,8,37,39,44            & 65  & 10          & 1,3,5,22,24,41,43,60,62,64 \\
46  & 8           & 1,3,5,22,24,41,43,45      & 66  & 8           & 1,3,9,32,34,57,63,65       \\
47  & 6           & 1,3,13,34,44,46           & 67  & 8           & 1,3,24,28,39,43,64,66      \\
48  & 6           & 1,10,21,27,38,47          & 68  & 7           & 1,3,13,34,55,65,67         \\
49  & 6           & 1,3,19,30,46,48           & 69  & 8           & 1,3,5,19,50,64,66,68       \\
50  & 7           & 1,3,14,25,36,47,49        & 70  & 8           & 1,3,26,30,40,44,67,69      \\
51  & 6           & 1,12,23,28,39,50          & 71  & 8           & 1,3,7,26,45,64,68,70       \\
52  & 6           & 2,10,13,39,42,50          & 72  & 8           & 1,6,8,35,37,64,66,71       \\
53  & 10          & 1,3,5,7,11,42,46,48,50,52 & 73  & 8           & 1,3,15,17,56,58,70,72      \\
54  & 6           & 1,10,24,30,44,53          & 74  & 8           & 1,7,13,30,44,61,67,73      \\
55  & 6           & 1,5,21,34,50,54           & 75  & 8           & 1,3,5,29,46,70,72,74       \\
56  & 8           & 1,3,7,26,30,49,53,55      & 76  & 8           & 1,3,14,18,58,62,73,75      \\
57  & 8           & 1,5,18,22,35,39,52,56     & 77  & 8           & 1,3,13,23,54,64,74,76      \\
58  & 8           & 1,3,7,26,32,51,55,57      & 78  & 6           & 1,12,17,61,66,77           \\
59  & 8           & 1,5,11,17,42,48,54,58     & 79  & 8           & 1,3,13,29,50,66,76,78      \\
60  & 7           & 1,3,19,30,41,57,59        & 80  & 8           & 1,3,13,34,46,67,77,79     
\end{tabular}

\bigskip

Small Values of $\psi_4 (n)$ 
\end{center}

Analogous to Question \ref{question 1}, one can ask if $\psi_{ \ell} (n) = o(n)$ 
for even $\ell \geq 4$.  We conclude the introduction by
noting that $\textup{rsat}(n , C_{ \ell +1} ) \leq n  \psi_{ \ell }(n)$
whenever these values exist (see 
Proposition \ref{subsum} or \cite{t}).


\subsection{Organization and Notation}

In Section \ref{key lemma section} we prove Theorem \ref{first theorem}.  The starting point will be a key lemma which is first 
applied in $\mathbb{Z}$ in Section \ref{subsection Z}.  Applying the key lemma in 
$\mathbb{Z}_n$ is done in Section \ref{subsection Zn}, which 
then culminates in a proof of Theorem \ref{first theorem}.
In Section \ref{graph section} we prove Theorem \ref{main theorem} and Corollary \ref{corollary}.  

Given even integers $m_1$ and $m_2$ with $m_1 < m_2$, let 
\[
[m_1 , m_2]_e := \{ m_1 , m_1 + 2 , m_1 + 4 , \dots , m_2 \}.
\]
Similarly, if $m_1 < m_2$ are both odd, 
\[
[m_1 , m_2]_o := \{ m_1 , m_1  + 2 , m_1 + 4 , \dots , m_2 \}.
\]
If $\Gamma$ is a group and $S \subseteq \Gamma$ is an 
inverse closed subset of $\Gamma$, then $\textup{Cay}( \Gamma, S)$ denotes the corresponding Cayley graph.  This is the graph with 
vertex set $\Gamma$.  Two distinct vertices $x$ and $y$ are adjacent if $xy^{-1} \in S$.  When $\Gamma$ is written 
additively, this last condition is $x -y \in S$.  In this paper 
$\Gamma$ will always be the cyclic group $\mathbb{Z}_n$.

\section{Key Lemma and Theorem \ref{first theorem}}\label{key lemma section}

\subsection{Statement and proof of Key Lemma}

In this section we prove a lemma that gives a formula
for a particular $\ell$-fold sumset.  
The initial set will be the union of three intervals of consecutive odd numbers, the last interval being a singleton.  We will then show that the corresponding $\ell$-fold sumset is 
an interval of consecutive even numbers, but with one interval of 
evens and one singleton removed (see equation
(\ref{eq:key lemma}) below).  

Let $\ell \geq 4$ be even, $\alpha \geq 1$, and $t \geq 1$ be integers.  
In $\mathbb{Z}$, define $I_1 = \{ 2 \gamma + 1 : 0 \leq \gamma \leq \alpha \}$, 
$I_2  = \{ 2 \alpha + 5 + 2 \theta : 0 \leq \theta \leq t \}$, $I_3 = \{ 2 \alpha + 4 t + 9 \}$, and $S^+ = I_1 \cup I_2 \cup I_3 $.   

\begin{lemma}[Key Lemma]\label{key lemma}
If $\alpha \geq 2t + 2 \ell - 2 $ and $M = \ell ( 2 \alpha + 4t + 9)$, then 
\begin{equation}\label{eq:key lemma}
\ell S^+ = [ \ell ,  M  ]_e \backslash ( [M - 2t - 2 , M  -2 ]_e \cup \{ M - 4t - 6 \} ).
\end{equation}
\end{lemma}
\begin{proof}
Given integers $0 \leq \beta_1 , \beta_2 , \beta_3 \leq \ell$ with $\beta_1 + \beta_2 + \beta_3 = \ell$, let
\[
J_{ \beta_1 , \beta_2 , \beta_3 } =  \beta_1 I_1 + \beta_2 I_2 + \beta_3 I_3.
\]
Observe that 
\[
\ell S^+ = \bigcup_{ \beta_1 , \beta_2 , \beta_3} J_{ \beta_1 , \beta_2 , \beta_3 }
\]
and our proof strategy will be to determine the union on the right.  
By definition,  
\begin{eqnarray*}
	J_{ \beta_1 , \beta_2 , \beta_3 } & = & 
	\{ \beta_1 + 2 ( \gamma_1 + \dots + \gamma_{ \beta_1 } ) + \beta_2 ( 2 \alpha + 5)  + 2 ( \theta_1 + \dots + \theta_{ \beta_2 } ) + \beta_3 ( 2 \alpha + 4t + 9)   \\
	& : & 0 \leq \gamma_i \leq \alpha , 0 \leq \theta_i \leq t \}  \\
	& = & \{ \beta_1 + \beta_2 ( 2 \alpha + 5) + \beta_3 ( 2 \alpha+ 4t + 9 ) + 2 \tau :
	0 \leq \tau \leq \beta_1 \alpha + \beta_2 t \}.   
\end{eqnarray*}
Using the equation $M =  \ell ( 2 \alpha + 4t + 9 )$, we can rewrite this expression as 
\begin{equation}\label{J eq1}
J_{ \beta_1 , \beta_2 , \beta_3 } = 
[
M - \beta_1 ( 2 \alpha + 4t + 8 ) - \beta_2 ( 4 t + 4) , 
M - \beta_1 ( 4t + 8 ) - \beta_2 ( 2t + 4) ]_e.
\end{equation}
We will use formula (\ref{J eq1}) often in the remainder of the proof.  
Strictly speaking, the $J_{ \beta_1 , \beta_2 , \beta_3}$'s 
are intervals of consecutive even numbers, but we 
will call them \emph{intervals} for short.  

Since 
\[
\min I_1 < \max I_1 < \min I_2 < \max I_2 < \min I_3 < \max I_3,
\]
we have 
\[
\min J_{ \beta_1 , \beta_2 , \beta_3 } < \min J_{ \beta_1' , \beta_2' , \beta_3' }
~~~~\mbox{and}~~~~
\max J_{ \beta_1 , \beta_2 , \beta_3 } < \max J_{ \beta_1' , \beta_2' , \beta_3' }
\]
whenever $(\beta_1 , \beta_2 , \beta_3 ) \prec ( \beta_1' , \beta_2' , \beta_3' )$.  Here $\prec$ is the lexicographic
ordering on 3-tuples with entries in $\{0,1, \dots , \ell \}$.  We 
extend this ordering to the $J_{ \beta_1 , \beta_2 , \beta_3}$'s
by saying $J_{ \beta_1 , \beta_2 , \beta_3 } \prec 
J_{ \beta_1 ' , \beta_2 ' , \beta_3 ' }$ if and only if
$( \beta_1 , \beta_2 , \beta_3 ) \prec ( \beta_1 ' , \beta_2 ' , \beta_3 ')$.  Now we consider these intervals from the 
largest, which is $J_{0,0, \ell } = \ell I_3 =   \{ M \}$, to the smallest, which 
is 
\[
J_{ \ell , 0 , 0 } = \ell I_1 =  [ M - \ell (2 \alpha + 4t  + 8 ) , M - \ell ( 4t + 8)]_e 
=
[ \ell , \ell (2 \alpha + 1) ]_e.
\]

The three largest intervals are $J_{ 0 , 0 , \ell } = \{ M \}$,  
\[
 J_{ 0 , 1, \ell - 1} = [ M - 4t - 4 , M  - 2t - 4 ]_e, 
~~~\mbox{and}~~~ 
J_{ 0 , 2 , \ell - 2 } = [ M - 8t - 8  , M - 4t - 8 ]_e.
\]
Thus,
\begin{equation*}\label{first three J}
J_{ 0 , 0 , \ell} \cup J_{0 , 1 , \ell - 1} \cup J_{0 , 2 , \ell - 2 } = 
[ M  - 8t - 8 , M]_e \backslash \left(  \{ M - 4t - 6 \} \cup [ M - 2t - 2 , M - 2 ]_e \right).
\end{equation*}
We will now show that the consecutive pairwise unions of the remaining $J_{ \beta_1 , \beta_2 , \beta_3}$'s form intervals of consecutive even numbers, 
some of which may overlap.  This implies that the entire union of the remaining $J_{ \beta_1 , \beta_2  , \beta_3}$'s is an interval of consecutive even numbers.  

Let $2 \leq \beta_2 \leq \ell - 1$ and consider 
$J_{ 0 , \beta_2 , \ell - \beta_2}$
and
$J_{ 0 , \beta_2 + 1 , \ell - \beta_2 - 1 }$.
The union of these two sets is the interval 
\[
[  \min J_{ 0 , \beta_2 + 1 , \ell - \beta_2 - 1 } , \max J_{ 0 , \beta_2 , \ell - \beta_2 } ]_e
=
[ M - ( \beta_2 + 1)( 4t + 4) , M - \beta_2 ( 2t + 4)]_e
\]
provided that 
\begin{equation}\label{ineq for J}
\max J_{ 0 , \beta_2 + 1 , \ell - \beta_2 - 1 } + 2 \geq \min J_{ 0 , \beta_2 , \ell - \beta_2}.
\end{equation}
This inequality is equivalent to 
\[
M - 2 t \beta_2 - 4 \beta_2 - 2t - 4 + 2 \geq M - 4 t \beta_2 - 4 \beta_2
\]
which can be rewritten as 
\begin{equation}\label{ineq 1}
2 t \beta_2 \geq 2t + 2.
\end{equation}
Since $\beta_2 \geq 2$ and $t \geq 1$, (\ref{ineq 1}) is true.  Thus, (\ref{ineq for J}) holds and  
\[
J_{ 0 , \beta_2 , \ell - \beta_2 } \cup J_{0 , \beta_2 + 1 , \ell - \beta_2 - 1} = [ M - ( \beta_2 + 1)( 4t + 4) , M - \beta_2 ( 2t + 4)]_e
\]
for $2 \leq \beta_2 \leq \ell - 1$.  We conclude that
\[
\bigcup_{ \beta_2 = 0}^{ \ell} 
J_{ 0 , \beta_2 , \ell - \beta_2 } 
=
[ M  - \ell ( 4t + 4) , M]_e \backslash \left(  \{ M - 4t - 6 \} \cup [ M - 2t - 2 , M - 2 ]_e \right)
\]
because $\min J_{ 0 , \ell , 0 } = M - \ell ( 4t + 4)$ by (\ref{J eq1}).    

We now complete the proof by showing that the following 
consecutive
pairwise unions are intervals of consecutive even numbers:
\begin{itemize}
	\item $J_{ \beta_1 , \ell - \beta_1 , 0 } \cup J_{ \beta_1 + 1 , 0 , \ell - \beta_1 - 1}$ for $0 \leq \beta_1 \leq \ell - 1$,
	\item $J_{ \beta_1 , \beta_2 , \ell - \beta_1 - \beta_2 } \cup J_{ \beta_1 , \beta_2 + 1 , \ell - \beta_1 - \beta_2 - 1 }$ for $0 \leq \beta_1 \leq \ell -1$ and $0 \leq \beta_2 \leq \ell - \beta_1 - 1$.
	\end{itemize}
Let $0 \leq \beta_1 \leq \ell - 1$.  The inequality 
$\max J_{ \beta_1+1 ,0 ,  \ell - \beta_1-1 }   +2 \geq  \min J_{ \beta_1 , \ell - \beta_1 , 0 } $
is equivalent to 
\[
\beta_1 ( \alpha - 2t - 2 ) + 2 \ell ( t + 1) \geq 2 t + 3.
\]
Since $\alpha - 2t - 2 \geq 0$ and $2 \ell ( t + 1) \geq 2t + 3$, 
this inequality holds.  Hence, 
$J_{ \beta_1 , \ell - \beta_1 , 0 } \cup J_{ \beta_1 + 1 , 0 , \ell - \beta_1 - 1}$ is an interval of consecutive even numbers.  

Next, let $1 \leq \beta_1 \leq \ell - 1$ and 
$0 \leq \beta_2 \leq \ell - \beta_1 -1 $.  The inequality 
\[
\max J_{ \beta_1 , \beta_2 + 1 , \ell - \beta_1 - \beta_2 - 1 }
+ 2 
\geq \min 
 J_{ \beta_1 , \beta_2 , \ell - \beta_1 - \beta_2 } 
\]
is equivalent to 
\[
2 \alpha \beta_1 + 2 \beta_2 t \geq 2 t+ 2.
\]
This holds since $\beta_1 \geq 1$ and $\alpha \geq 2t + 2$ which implies 
$2 \alpha \beta_1 \geq 2 ( 2t + 2)$.  Thus, 
\[
J_{ \beta_1 , \beta_2 , \ell - \beta_1 - \beta_2 } \cup
J_{ \beta_1 , \beta_2 + 1 , \ell - \beta_1 - \beta_2 - 1 }\]
is an interval of consecutive even numbers.

Since $J_{ \ell , 0 , 0 } = [ \ell , \ell ( 2 \alpha + 1)]_e$ and 
in particular $\min J_{ \ell , 0 , 0 } = \ell$, we have
\[
\ell S^+ = \bigcup J_{ \beta_1 , \beta_2 , \beta_3} = [ \ell , M] \backslash 
( \{ M - 4 t - 6 \} \cup [ M - 2t -2 , M - 2 ]_e ) 
\]
which completes the proof of Lemma \ref{key lemma}.  
\end{proof}


\subsection{Applying Key Lemma in $\mathbb{Z}$}\label{subsection Z} 

In this subsection we will apply Lemma \ref{key lemma}
to a certain choice of $I_1$, $I_2$, and $I_3$.  We begin by setting the parameters.  
Let $\ell \geq 4$ be even, $t \geq 1$, $\gamma$
be the unique integer that satisfies the inequality 
\[
- ( \gamma - 1) ( 2 \ell + 2) 
\geq \ell + 3 - 4 ( t + 2) 
\geq 
- \gamma ( 2 \ell + 2 ) + 1,
\]
$j = 2 \gamma + 1$, and $k \geq 4t + 2 \ell + 2 | \gamma |+2$.    
Define 
\begin{itemize}
	\item $I_1 = \{ 2 \theta + 1 : 0 \leq \theta \leq k + \gamma - 2t - 4 \}$,
	\item $I_2 = \{ 2k + j - 2 ( 2t + 4 ) + 2 \theta : 2 \leq \theta \leq t + 2 \}$, 
	\item $I_3 = \{ 2k + j \}$, and 
	\item $S^+ = I_1 \cup I_2 \cup I_3$.   
	\end{itemize}

We will apply Lemma \ref{key lemma} with $\alpha = k + \gamma - 2t - 4$ and $M = \ell ( 2k + j )$.  
The lower bound on $k$ implies that $\alpha \geq 2t + 2 \ell - 2$ so the hypothesis are satisfied.

With 
this choice of $\alpha$, we have $I_1 = \{ 2 \gamma + 1 : 0 \leq \gamma \leq \alpha \}$.  Furthermore, using the 
definition of $I_2$ given above, we have 
\begin{eqnarray*}
I_2 & = & \{ 2k + j  - 2 ( 2t + 4) + 2 \theta : 2 \leq \theta \leq t + 2 \} \\
& = & \{ 2 k + 2 \gamma - 4t - 8 + 5 + ( 2 \theta - 4) : 2 
\leq \theta \leq t + 2 \} \\
& = &  \{ 2 \alpha + 5 + 2 \theta : 0 \leq \theta \leq t \}.
\end{eqnarray*}
Likewise, 
\begin{eqnarray*}
I_3 & = & \{ 2k + j \} = \{ 2 ( k  + \gamma ) + 1 \}  \\
& = & \{ 2 ( k + \gamma - 2t - 4 ) + 4t + 9 \} = \{ 2\alpha + 4t + 9 \}.
\end{eqnarray*}

By Lemma \ref{key lemma}, 
\begin{equation}\label{eq:critical equation}
\ell S^+ = 
 [ \ell ,  M  ]_e \backslash
  ( [M - 2t - 2 , M  -2 ]_e \cup \{ M - 4t - 6 \} ) 
  \end{equation}
  where 
  \[
  M = \ell ( 2 \alpha + 4t + 9 ) = \ell ( 2k + j).   
  \]
Formula (\ref{eq:critical equation}) is critical to our proof of
Theorem \ref{first theorem}.  The next phase
will be to choose $n$ so that (\ref{eq:critical equation}) holds in $\mathbb{Z}_n$, and
$S^+$ is the ``positive half" of our complete symmetric $(\ell , 1)$-sum-free set in 
$\mathbb{Z}_n$.


\subsection{Applying Key Lemma in $\mathbb{Z}_n$ and the Proof of Theorem \ref{first theorem}}\label{subsection Zn} 	

In Section \ref{subsection Z} we applied 
Lemma \ref{key lemma} to a particular set $S^+$ viewed 
as a subset of $\mathbb{Z}$.  We will use $S^+$ 
to define a symmetric complete $( \ell , 1)$-sum-free set $S \subseteq \mathbb{Z}_n$. 

Assume that $\ell$, $t$, $\gamma$, $j$, and $k$
are chosen as in Section \ref{subsection Z}.    
Let 
\[
r = \ell + 3 - 4 ( t + 2) + \gamma ( 2 \ell + 2).
\]
The definition of $\gamma$ implies that
$r$ is the unique integer in $\{1,3,5 , \dots , 2 \ell  + 1 \}$ with 
\[
r \equiv \ell + 3 - 4 ( t + 2) ( \textup{mod}~ 2 \ell + 2 ).
\]
Lastly, let $n= ( 2 \ell + 2) k + r$.  Using the equations
defining $n$, $r$, and $j$, one can prove the following lemma by direct computation.  
The proof will be omitted.

\begin{lemma}\label{n equation}
	With $\ell, t , \gamma, j, k, r$, and $n$ 
	chosen as above, 
	\[
	\ell ( 2k + j ) = n - ( 2k + j - 2 ( 2 t + 4)  + 2).
	\]
	\end{lemma}

As before, we write $M = \ell ( 2k + j)$.  Let $I_1$, $I_2$, and 
$I_3$ be defined as in Section \ref{subsection Z}
and $S^+ = I_1 \cup I_2 \cup I_3$.  View $S^+$ as a subset 
of $\mathbb{Z}_n$ where, unless otherwise stated, 
we will always use the least residues $\{0,1, 2, \dots , n - 1 \}$.
Before proceeding further, we show
that $M < n$.  The reason that this is important is that because if
$M  < n$, then 
$\ell S^+$ is the same when the arithmetic 
is done in $\mathbb{Z}$ or $\mathbb{Z}_n$.  This allows us to apply the results of Section \ref{subsection Z}.  

By Lemma \ref{n equation}, 
\begin{eqnarray*}
n - M &  =&  2k + j - 2 (2t + 4) + 2  = 2k + j - 4t - 6 = 2 ( k + \gamma) - 4t - 5 \\
& \geq & 2 ( 4t + 2 \ell + 2 ) - 4t - 5= 4 t + 4 \ell - 1 .
\end{eqnarray*}
The right hand side is positive since $\ell \geq 4$ and so 
$M < n$.  Therefore, from (\ref{eq:critical equation}) in Section \ref{subsection Z}, 
\begin{equation}\label{l S plus}
\ell S^+ = 
 [ \ell ,  M  ]_e \backslash
  ( [M - 2t - 2 , M  -2 ]_e \cup \{ M - 4t - 6 \} ) 
  \end{equation}
  where $M =  \ell ( 2k + j)$ and now this is all in the abelian 
  group $\mathbb{Z}_n$. 
   Define 
  \[
  S^- = \{ n - s : s \in S^+ \} ~~
  \mbox{and}
  ~~
  S = S^+ \cup S^-.
  \]
  The set $S$ is a symmetric
subset of $\mathbb{Z}_n$ and we will now prove that 
\[
\ell S = \mathbb{Z}_n \backslash S.
\]
First we find the intervals that make up $S^-$.  
The set $S^-$ is the disjoint union of 
$n - I_1$, $n - I_2$, and $n - I_3$ where 
$n-I_i = \{ n - s : s \in I_i \}$.  
By Lemma \ref{n equation}, 
\[
n- I_3 = \{ n - ( 2k + j) \} = \{  \ell ( 2k + j ) - 2 ( 2t + 4) + 2 \} = \{ M - 4t - 6 \}.
\]
Similarly, 
\begin{eqnarray*}
n - I_2 & = & 
\{ n - ( 2k + j - 2 ( 2t + 4) + 2 \theta ) : 2 \leq \theta \leq t + 2 \}  \\
& = & 
\{ n - ( 2k + j  - 2 ( 2t + 4) + 2)   - 2 ( \theta - 1) : 2 \leq \theta \leq t + 2 \} 
\\
& = & \{ \ell (2k + j) - 2 \theta  : 1 \leq \theta \leq t + 1 \} 
\\
& = & [ M  - 2t - 2 , M - 2 ]_e.
\end{eqnarray*}
Lastly, 
\begin{eqnarray*}
n - I_1 & = & \{ n - 2 \theta - 1 : 0 \leq \theta \leq \alpha \} 
 =  [ n - 2 \alpha - 1 , n - 1]_e  \\
& = & [ n - 2 ( k + \gamma  -2 t - 4) - 1 , n -1]_e 
 =  [ M + 2 , n - 1]_e
\end{eqnarray*}
where in the last line we have used Lemma \ref{n equation}.  
Therefore, 
\begin{eqnarray*}
S^- & =& ( n - I_3) \cup (n - I_2) \cup ( n- I_1) \\
&=  & \{ M - 4t - 6 \}\cup [M - 2t - 2 , M - 2]_e \cup  [M+2,n-1]_e.
\end{eqnarray*}

Returning to (\ref{l S plus}) we now see that 
\begin{equation}\label{first lS plus}
\ell S^+ = [ \ell , M]_e \backslash ( ( n - I_2) \cup (n - I_3) ).
\end{equation}
Let $\mathcal{E} = \{0, 2,4, \dots , n - 1\} \subset 
\mathbb{Z}_n$.  Using 
(\ref{first lS plus}) and the equation $n - I_1 = [M+2,n-1]_e$, we can write $\mathcal{E}$ as the disjoint 
union 
\[
\mathcal{E} = [0 , \ell - 2 ]_e \cup \ell S^+ \cup (n - I_1) \cup ( n - I_2) \cup (n - I_3)
=
[ 0 , \ell - 2 ]_e \cup \ell S^+ \cup S^-.  
\]
By symmetry, if $\mathcal{O} = \{  1, 3, \dots , n - 2 \} \subset 
\mathbb{Z}_n$, 
then 
\[
\mathcal{O} = [ n - \ell + 2 , n-2]_o 
\cup \ell S^- \cup S^+.
\]
Putting these two together gives
\begin{equation}\label{first disjoint union}
\mathbb{Z}_n = [ 0 , \ell - 2]_e \cup S^+ \cup \ell S^- 
\cup \ell S^+ \cup S^- \cup [ n - \ell + 2 , n-2 ]_o.  
\end{equation}
We will now show that 
\begin{itemize}
    \item $[ 0 , \ell - 2]_e \cup \ell S^- \cup \ell S^+ \cup [ n - \ell + 2 , n-2 ]_o \subseteq \ell S$ and
    \item $S \cap \ell S = \emptyset$.  
    \end{itemize}

Since $S = S^+ \cup S^-$, we have 
$( \ell S^+ \cup \ell S^- ) \subseteq \ell S$.  
Now $S$ contains 1, 3, and $n - 1$.  Given 
$2z \in [ 0 , \ell - 2]_e$, we can write $2z$ as
\[
2z \equiv \underbrace{3 + \dots + 3}_{ \mbox{$z$ terms}} 
+ \underbrace{1 + \dots + 1}_{ \mbox{$l/2 - z$ terms}} 
+ \underbrace{(n-1) + \dots + (n-1)}_{ \mbox{$l/2$ terms}}
( \textup{mod}~n)
\]
and this sum is in $\ell S$.  Therefore, 
$[ 0 , \ell - 2]_e \subseteq \ell S$.  By 
symmetry,
$[n - \ell +2 , n-2 ]_o \subseteq \ell S$.  

Moving on to showing that $S \cap \ell S = \emptyset$, 
suppose for contradiction that $s \in S \cap \ell S$. 
By symmetry we may assume that $s \in S^+$.  
There are elements $s_i \in S$ such that 
\[
s \equiv s_1 + \dots + s_{ \ell} ( \textup{mod}~n).
\]
By relabeling we may assume $s_1 , \dots , s_r \in S^+$
and $s_{r+1} , \dots , s_{ \ell} \in S^-$ for 
some $0 \leq r \leq \ell$.  By definition of 
$S^-$,
there are elements $s_i' \in S^+$ such that 
$n - s_i' \equiv s_i ( \textup{mod}~n)$.  We can 
rewrite the above congruence as 
\begin{equation}\label{a congruence}
s \equiv
s_1 + \dots + s_r  -  ( s_{r+1}' + \dots + s_{\ell}' ) 
+ (\ell - r )n ( \textup{mod}~n)
\end{equation}
where $s_i , s_i' \in S^+$.  This is 
one point in our proof where we do not restrict to using 
the least residues $\{0,1, \dots , n - 1 \}$ for the elements of $\mathbb{Z}_n$.  The congruence (\ref{a congruence}) implies
\begin{equation}\label{integer equation}
s = s_1 + \dots + s_r - ( s_{r+1}' + \dots + s_{\ell}' ) 
+ \zeta n 
\end{equation}
for some integer $\zeta$.  Taking this integer equation modulo 2 gives 
\[
1 \equiv \ell + \zeta ( \textup{mod}~2).
\]
Recalling that $\ell$ is even, it must be the case 
that $\zeta$ is odd.  This, together with (\ref{integer equation}), implies 
\[
| s - s_1 - \dots - s_r +  s_{r+1}' + \dots + s_{\ell}'  | \geq n.
\]
This is a contradiction since the left hand side is at 
most $(2 \ell + 1)( 2k + j) < ( 2\ell + 2)k + r= n$ by the 
lower bound on $k$.  We have now shown that 
$S \cap \ell S = \emptyset$.  

Putting all of this together with (\ref{first disjoint union}) gives
\[
\mathbb{Z}_n = S \cup \ell S
\]
where the right hand side is a disjoint union.  
This proves Theorem \ref{first theorem}.


\section{Proof of Theorem \ref{main theorem} and 
Corollary \ref{corollary}}\label{graph section}

In this section we use Theorem \ref{first theorem} to 
prove that for any even integer $\ell \geq 4$,
there is an $n$-vertex $C_{\ell+1}$-saturated graph 
that is regular for all $n \geq 12 \ell^2 + 36 \ell + 24$.  Given 
$S \subseteq \Gamma$ where $\Gamma$ is an abelian group
and $\ell \geq 2$, let
\[
\mathcal{R}_{ \ell } (S)
= \{ s_1 + \dots + s_{ \ell} : s_i \in S 
\mbox{~and~} s_i + s_{i+1} + \dots + s_j \neq 0 
~\mbox{for all~} 1 \leq i < j \leq \ell  \}.
\]
We will combine Theorem \ref{first theorem} with 
the following proposition which is proved in \cite{t}.  

\begin{proposition}[\cite{t}]\label{subsum}
Let $\ell \geq 2$ be even and let $\Gamma$ be an 
abelian group with $| \Gamma | = n$.  
If there is a symmetric subset $S \subseteq \Gamma$ with 
\[
\mathcal{R}_{ \ell } (S) = \Gamma \backslash 
( S \cup \{ 0 \} ) ~~ \mbox{and} ~~
0 \notin ( \ell + 1) S,
\]
then the Cayley graph $\textup{Cay}( \Gamma , S)$ is
an $|S|$-regular $n$-vertex $C_{ \ell + 1}$-saturated 
graph.
\end{proposition}

For our application to graphs we need to consider 
$\mathcal{R}_{ \ell} (S)$ 
instead of $\ell S$.  The reason for this is that  
$C_{ \ell + 1}$-saturation requires at least 
one path of length $\ell$ between each pair of 
nonadjacent vertices.  The sums $s_1 + \dots + s_{ \ell }$ in 
$\ell S$ will be used to find these paths.  When 
a consecutive subsum is 0, we do not get a path of length $\ell$ 
(the path we are constructing returns to the initial vertex creating a cycle).  
Having made this remark, let us turn to the proof 
of Theorem \ref{main theorem}.  

\bigskip

\begin{proof}[Proof of Theorem \ref{main theorem}]
For even $n$, the graph $K_{n/2,n/2}$ is $C_{ \ell + 1}$-saturated whenever $n \geq \lceil \frac{ \ell + 1}{2} \rceil$.  
Let $n$ be an odd integer such that $n \geq 12 \ell^2 + 36 \ell + 24$.
Suppose that $n \equiv r ( \textup{mod}~ 2 \ell  + 2 )$ 
where $r \in \{1, 3, \dots , 2 \ell + 1 \}$ in $\mathbb{Z}_n$.

\begin{claim}\label{t claim}
There is an integer $t \in \{1,2, \dots ,  \ell  + 1 \}$ such that 
\[
\ell + 3 - 4 ( t + 2) \equiv r ( \textup{mod}~ 2 \ell  + 2).
\]
\end{claim}
\begin{proof}[Proof of Claim \ref{t claim}]
Suppose that $\ell + 3 - 4 (t_1 +2) \equiv \ell + 3 - 4 (t_2 +1) ( \textup{mod}~2 \ell + 2 )$ where $t_1,t_2 \in \{1,2, \dots , \ell + 1 \}$.
This congruence is equivalent to $2 t_1 \equiv 2 t_2 ( \textup{mod}~ \ell + 1)$.  Since $\ell +1 $ is odd, we have 
$t_1 \equiv t_2 ( \textup{mod}~ \ell + 1)$.  Since 
$t_1 ,t_2$ are integers in $\{1,2, \dots , 
\ell +1 \}$, it must be the case that $t_1 = t_2$.
Thus, 
\[
\{ \ell + 3 - 4 ( t + 2) ( \textup{mod}~ 2 \ell + 2) 
: t \in \{1,2, \dots , \ell + 1 \} \} =
\{1 , 3, \dots , 2\ell + 1 \}.
\]
\end{proof}

\bigskip

By Claim \ref{t claim} we can choose an integer 
$t \in \{1,2, \dots , \ell  + 1 \}$ such that 
\[
\ell + 3 - 4 ( t + 2) \equiv r  ( \textup{mod}~2 \ell + 2).
\]
We choose $\gamma$ so that 
$r = \ell + 3 - 4 (t + 2) + \gamma ( 2 \ell + 2 )$
where $\gamma$ is the unique integer satisfying 
\[
- ( \gamma - 1) ( 2\ell + 2) \geq \ell + 3 
- 4 (t + 2 ) \geq - \gamma ( 2 \ell + 2) + 1.
\]
From the definition of $\gamma$, 
$r$ is the unique integer in $\{1,3, \dots , 2 \ell + 1 \}$
with $n \equiv r ( \textup{mod}~ 2 \ell + 2 )$.  
The lower bound on $n$ implies that 
there is an integer $k \geq 6 \ell + 12$ such that 
$n = (2 \ell + 2 ) k + r$.  This $k$ will 
satisfy the inequality $k \geq 4 t + 2 \ell + 2 | \gamma |+2$ which 
is needed in order to apply Theorem \ref{first theorem}.
We take a brief moment to check that this is indeed the case.  

From the definition of $t$, we know that $t \leq \ell + 1$. 
The assumption on $\gamma$ implies that 
\begin{equation}\label{gamma ineq}
\gamma \leq \frac{ 4 ( t + 2) - \ell - 3 }{ 2 \ell + 2 } + 1 
\leq 
\frac{ 5 \ell + 1}{ 2 \ell + 2} < \frac{ 8 \ell }{ 2\ell } = 4
\end{equation}
where we have used the fact that $t \leq \ell + 1$ and $\ell \geq 4$.  Hence, 
\[
4t + 2\ell + 2 | \gamma| +2
\leq 4 ( \ell + 1) + 2 \ell + 6+2 = 6 \ell + 12
\]
and so since $n \geq 12 \ell^2 + 36 \ell + 24$, we have
$k \geq 6 \ell + 12 \geq 4 t + 2 \ell + 2 | \gamma |+2$.

Taking $j = 2 \gamma + 1$, we apply Theorem \ref{first theorem}
to obtain a symmetric subset $S \subseteq \mathbb{Z}_n$ for 
which $\mathbb{Z}_n$ is the disjoint union of $S$ and 
$\ell S$.  To apply Proposition \ref{subsum}, we need to prove that 
for this $S$, 
\begin{equation}\label{R equal}
\mathcal{R}_{ \ell} (S) = \mathbb{Z}_n \backslash ( S \cup \{ 0 \} ) 
\end{equation}
and $0 \notin ( \ell + 1) S$.  For the latter, if 
$0 \in ( \ell + 1) S$, then there are elements 
$s_1 , \dots, s_{ \ell + 1}$ such that 
$s_1 + \dots + s_{ \ell + 1} \equiv 0 ( \textup{mod}~ n)$.  Because $S$ is symmetric, $- s_{ \ell + 1} \in S$ and 
therefore, this congruence implies $- s_{ \ell + 1} \in \ell S \cap S$, a contradiction.  

We will complete the proof of Theorem \ref{main theorem}
by showing that (\ref{R equal}) holds.  
To prove this, it is enough to show that every 
element $x \in\ell S \backslash \{ 0 \}$ can be written as 
\[
x \equiv s_1 + \dots + s_{ \ell}  ( \textup{mod}~n), ~s_i \in S
\]
where $s_i + \dots + s_j \not \equiv 0  ( \textup{mod}~n)$ 
for all $1 \leq i < j \leq \ell$.  Critical to this is the structure of the sets $S^+$ and $S^-$ in the proof of Theorem \ref{first theorem}.  From (\ref{first disjoint union}) we know
\begin{equation}\label{a disjoint union}
\mathbb{Z}_n = 
\left(
[ 2 ,  \ell - 2 ]_e \cup \ell S^- \cup \ell S^+ 
\cup [ n - \ell + 2 , n - 2 ]_0 
\right)
\cup ( S \cup \{ 0 \} ).
\end{equation}
A sum $s_1 + \dots + s_{ \ell }$ in $\ell S^+$ is 
also a sum in $\mathcal{R}_{ \ell }(S)$ because
$S^+ \subset \{1,3, \dots , \ell ( 2k + j ) \}$ and 
$\ell ( 2k + j ) < n$ so that 
all consecutive subsums of sums in $\ell S^+$ 
are contained in
\[
\{1, 2, \dots , \ell ( 2k + j ) \}.
\]
By symmetry there is no consecutive 
subsum of $s_1 + \dots + s_{ \ell } \in \ell S^-$ 
that is 0 modulo $n$. 
Thus, $( \ell S^+ \cup \ell S^- ) \subseteq \mathcal{R}_{ \ell } (S)$.  It remains to express
each element in $[2 , \ell -2 ]_e \cup [ n - \ell + 2 , n - 2]_o$
as a sum in $\mathcal{R}_{ \ell }(S)$.  If
$2 u \in [ 2 , \ell - 2]_e$,
then 
\[
2u \equiv 
\underbrace{- ( 2k+j) - \dots - ( 2 k+j) }_{  \frac{\ell -2u}{2}~\mbox{times}}
+ \underbrace{1 + \dots + 1}_{2u~\mbox{times}}
+ \underbrace{(2k +j) + \dots + ( 2k + j) }_{ \frac{\ell -2u}{2}~\mbox{times}} ( \textup{mod}~n).
\]
The first $\frac{ \ell - 2u }{2}$ terms are in $S^-$ and the remaining terms are in $S^+$.  Suppose, for contradiction, that some consecutive subsum is 0.  Any such sum must use 
terms from both $S^-$ and $S^+$ since $\ell ( 2k + j)  < n$.  
Assume that 
\[ 
- \theta_1 ( 2k + j) + 2u + \theta_2 ( 2k + j) \equiv 0 
( \textup{mod}~n)
\]
where $1 \leq \theta_1 , \theta_2 \leq \frac{ \ell - 2u }{2}$. 
If $\theta_2 \geq \theta_1$, then 
\[
2 u + ( \theta_2 - \theta_1 ) ( 2k + j) \equiv 0 ( \textup{mod}~n).\]
This is a contradiction because the left hand side 
is a least residue in the set 
\[
\{ 2 , 4, \dots , \ell ( 2k + j) \}
\]
and $\ell ( 2k  +j ) < n$.  
When $\theta_1 > \theta_2$, we have 
\[
 ( \theta_1 - \theta_2) ( 2k + j) \equiv 2u ( \textup{mod}~n).
 \]
 This is a contradiction as before because now the 
 left hand side is a least residue in 
 $\{ 2k +j ,  2 ( 2k + j ) , \dots , \ell ( 2k + j) \}$ while
 the right hand side is in $\{ 2 , 4, \dots , \ell - 2 \}$,
 and \[
 \ell - 2 < 2k + j < \ell (2k + j ) <n.
 \]
  The conclusion is that $[ 2 , \ell - 2]_e  \subseteq \mathcal{R}_{ \ell } (S)$ and by symmetry, $[ n - \ell + 2 , n - 2]_0$ is also 
 a  subset of $\mathcal{R}_{ \ell } (S)$.  Recalling 
 that 
 $ (S^+ \cup S^-) \subseteq \mathcal{R}_{ \ell }(S) \subseteq \ell S$
 and (\ref{a disjoint union}), we see that $\ell S \cap S = \emptyset$ implies 
 \[
 \mathcal{R}_{ \ell} (S) = \mathbb{Z}_n \backslash ( S \cup \{ 0 \} ).  \]
 By Proposition \ref{subsum}, $\textup{Cay}( \mathbb{Z}_n , S )$ is a $C_{ \ell + 1}$-saturated graph with $n$ vertices 
 and is $|S|$-regular.  
 \end{proof}

\bigskip

We end this section by demonstrating the upper bound 
on $\textup{rsat}(n, C_{ \ell + 1})$ that our construction 
gives.  
Recalling that $n$ is odd and $n = ( 2 \ell + 2) k + r$,
a short computation gives $|S| = 2 ( k + \gamma - t - 1 )$.  
The inequalities $\gamma \leq 3$ (see (\ref{gamma ineq})) and $t \geq 1$ imply that
\[
|S|  = 2 ( k + \gamma - t - 1 ) \leq 2 ( k + 1).
\]
Using the equation $k = \frac{n - r }{ 2( \ell + 1) }$ and 
this upper bound on $|S|$, we have for any even $\ell \geq 4$ 
and odd $n \geq 12 \ell^2 + 36 \ell + 24$,
\[
\textup{rsat}(n , C_{ \ell + 1} ) \leq \frac{n^2}{ 2 ( \ell + 1) } + n. 
\]


\end{document}